Routine Controversies:

Mathematical Challenges in Mersenne's Correspondence

Controversies have become in the 1980s a privileged topic in history of science.[1] This interest goes with the demythification—some might say the disenchantment—of science, or perhaps more precisely of its traditional public values (integrity, commonality, free exchange of ideas, absence of personal and political interests, objectivity) in the construction of concepts and in the statement of results. For people, scientists and historians alike, who believed that these values realistically described science, controversies were epiphenomena, insignificant testimonies of local malfunctioning, human, petty, issues with no important bearing on valuable scientific practice. More recently, on the contrary, controversies seem to have become the places to look at to capture the authentic structure of truth-making processes: "Controversies are the engine of intellectual progress in philosophy, science, theology, the arts, and other domains,"[2] claims the website of the International Association for the Study of Controversies. This position might accompany an otherwise positivistic view of science and scientific activities—polemics being seen as the tool towards a progress the existence or features of which remain unquestioned. Or it can be used to undermine oversimple ideas on progress, truth, scientific facts and reality. Steven Shapin and Simon Shaffer wrote for instance:

"Historical instances of controversy over natural phenomena or intellectual practices have two advantages, from our point of view. One is that they often involve disagreements over the reality of entities or propriety of practices whose existence or value are subsequently taken to be unproblematic or settled. […] Another advantage afforded by studying controversy is that historical actors frequently play a role analogous to that of [the stranger in anthropology]: in the course of controversy they attempt to deconstruct the taken for granted quality of their antagonists's preferred beliefs and practices, and they do this by trying to display the artifactual and conventional status of those beliefs and practices."[3]

Controversies would thus help to reveal actors and society otherwise neatly dissimulated in the final display of facts or theorems.

Topic in themselves or historiographical tool, controversies, however, are usually not historicized. Punctual studies pointing out to some of the available models scientific disputes have borrowed—in particular duels, or advertisements in a competitive market—suggest decisive interactions between the social and intellectual aspects of controversies. However, controversies are usually associated with disruption or changes of paradigms, and rarely integrated in the social order of so-called normal science. While born as a tool for the historian in the steps of the strong program in the sociology of science, they have been little submitted to its reflexivity constraint: is this tool itself dependent of the disciplines, of the times, of the places? How is it related to the cognitive structure of the practices under study?

[1] For significant examples of this new focus on controversies, see: Martin J. S. Rudwick, *The Great Devonian Controversy: The Shaping of Scientific Knowledge among Gentlemanly Specialists* (Chicago:University of Chicago Press, 1985); Steven Shapin & Simon Shaffer, *Leviathan and the Air-Pump* (Princeton: Princeton University Press, 1985); James Secord, *Controversy in Victorian geology : the Cambrian-Silurian dispute* (Princeton: Princeton University Press, 1986); Bruno Latour, Pasteur et Pouchet: hétérogénèse de l'histoire des sciences, in Michel Serres (dir.), *Éléments d'histoire des sciences* (Paris:Bordas, 1989), 422-455.

[2] http://www.tau.ac.il/humanities/philos/iasc/, consulté le 19 avril 2011.

[3] Shapin & Shaffer, *op. cit. in* n.1, 7.

My aim here is to describe controversies as part and parcel of a correspondence network, loosely defined as that of Marin Mersenne (a detailed definition of my corpus will be given later). This correspondence will be here considered as an institution[4]—one among the first scientific organizations aiming at cooperative work according to Baconian ideals, but also one in which mathematical challenges, many of them leading to controversies or expressing them, were the rule, not the exception: as we shall see, these challenges operated here both as mise en scène of methodological oppositions—not as their spontaneous outcome—, and as links in a mathematical environment structured in letters around the resolution of problems.

## I—Marin Mersenne's correspondence as an institution

Marin Mersenne's scientific activities can be traced back to the first decade of the seventeenth century, but the correspondence that has survived (and has been published in 17 volumes[5]) displays a significative increase in the middle of 1630s. This period also corresponds to the establishment of a Parisian academy "entirely devoted to mathematics,"[6] gathering in the fifteen years or so before Mersenne's death such famous participants as the Pascals (father and son), Gilles Personne de Roberval, Pierre Gassendi, Bernard Frenicle de Bessy or Claude Mydorge. But Mersenne dreamt of the correspondence itself as a kind of academy:

> "Je voudrais que nous eussions une telle paix que l'on put dresser une Académie, non dans une seule ville comme on fait icy et ailleurs, mais sinon de toute l'Europe, du moins de toute la France, laquelle entretiendrait ses communications par lettres, qui seroient souvent plus profitables que les entreparlers, ou l'on s'echauffe souvent trop à contester les opinions que l'on propose, ce qui en destourne plusieurs."[7]

[4] Here, an institution designates "any form of social organization that links values, norms, models of relations and behaviours, roles,," following a description of G. Balandier, quoted and commented by Jacques Revel, L'institution et le social, in Bernard Lepetit (dir.), *Les formes de l'expérience* (Paris: Albin Michel, 1995), 63-84.

[5] Marin Mersenne, *Correspondance*, éd. Marie Tannery, René Pintard, Cornelis de Waard, Bernard Rochot, 17 vols (vol. I, Paris: Beauchesne, 1932; II-III, Paris: PUF, 1945-1947; vol. III, 2e éd. et vols IV-XVII, Paris: CNRS, 1955-1988).

[6] Mersenne, *op. cit. in* n.5, V, 209.

[7] Mersenne, *op. cit. in* n.5, V, 301-302.

We see then that our very topic—the role of controversies in the functioning of science—was thematized by Mersenne (as it was of course by other theoreticians of scientific life as Francis Bacon), his wish being to engineer his correspondence as a social solution to the problems they raise. This solution, putting negative emotions aside by the mediation of writing, was apparently inefficient, if we keep in mind the multiple and bitter controversies of these decades, for example between Pierre Fermat and René Descartes or between Jean de Beaugrand and Girard Desargues, to mention two well-known cases. To study the functioning of the correspondence as an institution of science in itself, and to understand the place of controversies in it, is thus specially intriguing.

For sake of place, I shall limit myself here to the part of correspondence dealing with questions on numbers. The traditional classification of mathematics that opposes geometry to arithmetic, the first being devoted to magnitudes and figures, the second to numbers, was indeed transformed with the emergence of new areas like algebra or mechanics, but is still active in Mersenne's correspondence. For instance, Mersenne wrote to Christiaan Huygens: "Ces petites gentillesses de nombres sont trop scabreuses pour vous en entretenir, la Géométrie est plus joyeuse,"[8] or Fermat to Mersenne: "J'ai trouvé beaucoup de sortes d'analyses pour divers problèmes, tant numériques que géométriques."[9] About 160 letters then, written by or to forty-six different persons, devote at least a few lines, sometimes a dozen of pages, to questions on numbers.

Letters are both social links and texts. It is noticeable that neither aspect is here obvious, nor are they independent. Some letters sent to Mersenne are more or less explicitly addressed to somebody else, complexifying or even destroying the expected interpretation of a letter as a simple link between the sender and the addressee. In a letter to Mersenne, for instance, Fermat shortens an explanation to the point of obscurity, adding: "J'en dis assez pour me faire entendre de M. Frenicle;"[10] a link between the sender, Fermat, and the real addressee in the shadow, Frenicle, is thus indicated textually by refering to a tacit knowledge that the official addressee does not share. In some cases, only part of an original letter has survived, or it has been rewritten with different notations, and this material change testifies a subtle activity of textual appropriation by specific subgroups. Moreover, several letters are drafts of manuscripts or work, with no other copy existing:

[8] Mersenne, *op. cit. in* n.5, XV, 30. For a more precise description of the topics involved here, see Catherine Goldstein, L'arithmétique de Fermat dans le contexte de la correspondance de Mersenne: une approche micro-sociale, *Annales de la Faculté des sciences de Toulouse* XVIII, n.spécial, 2009, 25-57. However, key features, in particular the role of problems in challenges, would also be relevant for geometry.

[9] Pierre Fermat, *Œuvres*, éd. Paul Tannery, Charles Henry, 4 vols+supplément, éd. Cornelius de Waard (Paris: Gauthier-Villars, 1891-1922), II, 5.

[10] Fermat, *op. cit. in* n.9, II, 216.

they are thus not only registering mathematics that would be produced or transmitted in more recognizable forms elsewhere; they are defining the production and transmission of mathematics.[11] All these aspects are historiographically significant: we cannot study the correspondence as a social network while ignoring what is actually written in the letters; nor can we use them as a pure source of information, on mathematical results or debates, without taking into account their paratextual and social features.[12] The problem is particularly acute in our study, because controversies may have been exacerbated, or on the contrary smoothed, even annihilated (as Mersenne wished) by the medium itself.

Acknowledging these issues is decisive for the legitimacy of our project itself, that is the possibility of studying these letters as an institution. An obvious metaphor suggests to see it as a network linking persons and to use then network analysis with various characteristics attached to these persons (geography, social position, skills, …) or to their links (degree of intimacy, …) to describe it. Such an approach would precisely make us err: the first impression is one of chaos, with a cacophony of voices and situations, up to the point to break apart the unicity of the whole.[13] According to Mersenne's wishes, indeed, correspondents come from a variety of places and countries, offering a variety of mathematical environments, from regular and well informed meetings to almost complete isolation. The relation to the dynamic of work and communication is never simple. Paris may appear as the place to be or to avoid, local academies may fire controversies or help to flatten them: in the controversy about tangents between Fermat and Descartes, both have their Parisian group of champions replaying and refueling the arguments. Social positions go from the highest levels of the aristocracy (Elisabeth of Bohemia or Louis de Valois) to modest monks and craftsmen, with again multifarious effects.[14] Being to high in status might hamper the communication of the most delicate scientific points as much as being to low. Status outside mathematics may be used to reinforce a scientific information, scientific talent to reinforce a social status, and both can interfere to the detriment of the other: "J'ay si peu de loisir,

---

[11] The impact of this issue on our knowledge of Fermat, for instance, is discussed in Goldstein, *op. cit. in* n.8.

[12] These issues have already been partially tackled, see by David Lux & Harold Cook, Closed Circles or Open Networks? Communicating at a Distance during the Scientific Revolution, *History of Sciences*, xxxvi (1998), 179–211 and, on the need to study together the scientific production and the communication structures, Jeanne Peiffer, Faire des mathématiques par lettres, *Revue d'histoire des mathématiques* 4(1998), 143-157. But I also want to emphasize here the mathematical production *as* part of the communication structure.

[13] See Aurélien Ruellet, *Le réseau des correspondants de Marin Mersenne*, Mémoire de Master, Université de Tours, 2004-2005.

[14] See Catherine Goldstein, L'honneur de l'esprit: de la République des mathématiques, in *Dire et vivre l'ordre social en France sous l'Ancien Régime*, textes réunis par Fanny Cosandey, Paris: Ed. de l'EHESS, 2005, p. 191-230.

complaints de Beaune, que cela me feroit volontiers souhaiter une autre condition que la miene [sic], si les affaires domestiques le pouvoient permettre."[15]
Commitment to mathematics, competencies and knowledge again describe the whole spectrum. Some persons, like Descartes, seem to have almost all answers, some others, like Thibaut, none. And once more, a lot of graduations exist: a Bergerac doctor, Theodore Deschamps, can speak to Mersenne with relevant authority about the difference between the properties of numbers which depend of the decimal writing and those who do not, or muse that, "du temps où [il était] escolier et où les mathématiques étaient [sa] plus grande passion,"[16] he also imagined to use lines in algebra as he discovered Anderson, Descartes and Billy have done. Claude Bredeau, an important interlocutor on musical questions, has nonetheless difficulties to read Mersenne's La vérité des sciences,

"pour n'avoir pas esté enseigné aux Mathématiques; pour auxquelles entendre quelque chose j'ay tousjours cru qu'il fallait sçavoir les Elements d'Euclide par coeur. Je n'ay eu aucun qui m'ayt tracé le chemin, temps ni loisir pour ce faire."[17]

Renaissance's editions of Archimedes, Apollonius, Diophantus and all their modern commentaries are courant knowledge for Fermat, while Frenicle does not know or does not use symbolic algebra. Terminology is not shared, meanings of some questions are subject to misunderstandings, originality of a procedure is put into doubt because of references to a different set of sources. For instance Sainte-Croix proposed in 1638 a problem about a "trigone tétragone," which is interpretated differently by various members of the network.[18]

As Lisa Sarasohn rightly put out,[19] we are confronted to a complex redistribution and exchanges and too quick or too simplistic models are condemned to fail. I already argued that the whole procedure is not describable as a transmutation of intellectual powers into wordly ressources and vice-versa. Neither is it a pure distributive system of finite products, a market of mathematical facts. Taken in isolation, certain documents certainly suggest one or the other of these analyses, but the intricacy of relationships is not dissolved through a clear-cut hierarchical system or a perfect specialization of tasks. While Fermat possesses most scientific answers and a quite comfortable position at the Parlement in Toulouse, and thus could seem to draw almost no advantage of being in the network, he expresses eagerly his wish to participate and the stimulation he found there.
It is then legitimate to ask whether a unique social space is concerned or several loosely connected or even if each pair of correspondents constitute a small private world of its own.

---

[15] René Descartes, *Œuvres*, ed. C. Adam et P. Tannery, Paris:Vrin, III, 537.
[16] Mersenne, *op. cit. in* n.5, XIII, 35-36.
[17] Mersenne, *op. cit. in* n.5, I, 496.
[18] Descartes, *op. cit. in* n.15, II, 158-161.
[19] Lisa Sarasohn, Nicolas-Claude Fabri de Peiresc and the Patronage of the New Science in the 17th Century, *Isis* 84, 1993, p. 70-90.

The picture changes when mathematics is taken into account—as links between persons or letters, and also as nodes. Problems on aliquot parts—that is on the computation of the divisors of integral numbers, their sums and properties—are alluded to for instance in more than fifty letters and half of the correspondents discuss them. They cement the group to the point of providing it with metaphors in various circumstances and even with private jokes: Deschamps uses the analogy with aliquot parts to discuss the infiniteness of "what is from itself," Mersenne compares the fecondity of François de la Noue, an important member of the Minim order to whom he dedicated his 1644 his Cosmographia astronomica to that of a number with many aliquot parts. Brief apparitions of certain questions, for instance the first case of Fermat so-called Last Theorem, in letters mainly devoted to different matters, reinforce the coherence of the correspondence. In this respect, we are thus confronted to a mathematical institution (in the sense of an institution made of mathematics, even more than of mathematicians), which had to solve the disruptive questions of the diversity we just explain. This is why we need to understand more closely the characteristics of the mathematics done in the correspondence and the way it circulates and transforms.

II—The organization of mathematical work

The first important feature is the shape of the statements. Let us look at a few typical examples: `Que [Fermat] vous envoie un nombre parfait qui ait 20 lettres ou le prochainement suivant" asks Frenicle to Mersenne; or Fermat to his former colleague Pierre de Carcavi: "Pour exciter par mon exemple les savants du pays où vous êtes, je leur propose de trouver autant de triangles rectangles en nombres qu'on voudra de même aire" or again:

"La seconde question [de M. de Saint-Martin] est celle-ci: un nombre étant donné, déterminer combien de fois il est la différence des côtés d'un triangle qui ait un quarré pour différence de son petit côté au deux autres côtés. Le nombre qu'il donne est 1 803 601 800. Je réponds qu'en l'exemple proposé il y a 243 triangles qui satisfont à la question et qu'il ne peut y en avoir davantage."[20]

The emphasis thus is on problems. I refer here to the standard distinction between problems and theorems, again inherited from Antiquity, in particular from Proclus's commentary on Euclid's Elements: in his 1691 Dictionnaire de mathématiques, Jacques Ozanam describes that « le Theoreme est une proposition speculative qui exprime les proprietez d'une chose, […] le Probleme est une proposition qui tend à la pratique. » Problems are at the core of the exchanges: Mersenne requires them from his correspondents, and we have traces of participants reformulating results under the form of a problem to satisfy him and the rules of the communication; bristling the problem with superficial complications to hide its main point to ignorant people is also two a penny.[21] Theorems when mentioned are facts to be used, eventually to be commented, more than results of a necessary demonstrative procedure.

[20] Fermat, *op. cit. in* n.9, II, 185, 248-249 and 250 respectively.

Answer to a problem is an explicit construction, in arithmetical questions a number. Numbers, particularly if they are big, that is not immediately accessible through immediate inspection, are all by themselves proofs: proofs of knowledge, proofs of success, both and together at the scientific and at the social level. Finding them, indeed, is the difficult part, verification (proof in the current sense) being normally the easiest part:

"Vous pouvez vous assurer que j'en possède absolument la méthode, […] et pour vous montrer jusques où va la connaissance que j'en ai, le quarré de 8, qui est 64, se peut disposer en autant de façons différentes [en carré magique] et non plus qu'il y a d'unités en ce nombre, 1004 144 995 344."[22]

This structure is supported by a general stratification of work: problems and their answers are bound to the application of a method, which allows the mathematician, from some principles (given or to be found and which are the typical place where theorems appear) to derive a rule (or a set of rules) through which to gain access to the concrete wanted solution. The method is supposed to be applicable to a certain category of problems, as large as possible of course, but testing it on concrete problems is at the end the only admitted means of deciding its value. Rules may even appear as a pis-aller, as for instance Descartes answers to Sainte-Croix' questions by a rule which is a non-practical one, because actual numbers would be too difficult to obtain concretely.[23] The chore of mathematical activities here is not to delimitate a certain corpus to be explored completely and established by fixed argumentative deductions, but to possess a method to solve any problem (or at least a range of them) which could be asked. The subjacent idea, expressed by some of the participants, is that you cannot solve a priori every problem in the mathematical world, because they are infinite in numbers, but we can and have to derive from a systematic procedure, a guide for your own practice, rules and solutions to be asked. Even in a closer proximity, this organization prints its mark. For instance, Mersenne writes to Saint-Martin around 1640:

---

[21] Catherine Goldstein, L'expérience des nombres de Bernard Frenicle de Bessy, *Revue de synthèse* (IV) 122 (2001), 425-454. As said before, the situation is not different for geometry, see for instance the problem sent by Florimond de Beaune circa 1639, which is: to find a curve such that the ratio of the subtangent to the ordinate is equal to the ratio of a constant line to the ordinate of the bissectrix of the angle of the axes decreased by the ordinate of the curve.

[22] Fermat, *op. cit. in* n.9, II, 195.

[23] Descartes, *op. cit. in* n.15, II, 162: "pour ce que ne sçay pas combien longs pourront estre les premiers [nombres solutions] qu'on rencontrera, j'ayme mieux mettre icy une regle […] que de m'arester moy mesme a faire le calcul qui est necessaire pour les chercher."

"Je vous prie que nous nous remettions un peu pour trouver combien un nombre proposé a de parties aliquotes, par exemple combien en a 49 000; et aussi quelle est la somme des parties aliquotes sans qu'il soit besoin de les compter. […] Je désire que vous m'apreniez votre méthode certaine."

To this, Saint-Martin only answers back: "Je trouve que votre nombre 49000 a 335 parties aliquotes et la somme de ses parties est 3255."[24] The method appears as a personal trade-mark of subtlety and talent, the rule a step by step procedure which would correspond at a possible level of intercommunication, the explicit solution as the minimal answer-unit, both at the very end, also as what is sought-for.

A typical example, which witnesses the coordination of several themes, including the time-factor one, is given in a the letter from Descartes to Mersenne:[25]

"Pour la façon dont je me sers à trouver les parties aliquotes, je vous diray que ce n'est autre chose que mon Analise, laquelle j'applique à ce genre de questions, ainsi qu'aux autres, & il me faudrait du temps pour l'expliquer en forme d'une règle, qui pust estre entenduë par ceux qui usent d'une autre méthode. Mais j'ay pensé que, si je mettois icy une demi-douzaine de nombres, dont les parties aliquotes fissent le triple, vous n'en feriez peut-estre pas moins d'estat, que si je vous envoyois une regle pour les trouver. C'est pourquoi je les ay cherchez, & les voicy:
30 240, dont les parties font … 90720
32760, dont les parties font … 98280
23569920, dont les parties font …70709760 […]."

Even later in the century, we still find traces of this structure. About the problem to solve in integers $NX^2+1=Y^2$, N fixed, John Wallis writes:[26]

"J'ai jugé bon de taire également les méthodes, soit de vous, soit de moi pour obtenir par induction le premier carré et sa racine […] je n'ai guère vu de moyen d'exposer clairement ces méthodes en sorte qu'elles soient facilement comprises par autrui, sans un appareil de mots et d'exemples. Quant au centre de gravité […] vous verrez qu'il manque beaucoup de choses que je vous ai déjà exposées là-dessus. […] J'ai préféré énoncer tout d'abord ces questions à Fermat sous forme de problèmes."

This type of exchange fits well in the framework of mutual obligations and gifts among extremely different participants described above.[27] It is also very coherent with one type of output of the correspondence, the books written by Mersenne that are made of a juxtaposition of such variegated pieces, with few attempts to arrange them in a hierarchy or even to develop argumentative lines to eliminate some of the answers. However, all of by itself, it does not seem to provide much stimulation to search for new results or arguments. When asked in 1631, to provide numbers besides 120 which are "double numbers" (that is, equal to half the sum of their divisors), Descartes answered that he has no idea and passes on. The way to force attention and work on a problem is to present it as a challenge—in particular to propose problems that are not open questions, but for which an answer, or the answer, is known by some participants.

[24] Fermat, *op. cit. in* n.9, IV, 69-70.

[25] Mersenne, *op. cit. in* n.5, VII, 345= Descartes, *op. cit. in* n.15, II, 250-251.

[26] John Wallis to William Brouncker, *Commercium epistolicum* (trad. P. Tannery), in Fermat, *op. cit. in* n. 9, III, 426.

[27] There is of course a striking analogy with what has been described as the new regime of experimental truth, with emphasis on individual instances and the possibility that this opens to tolerant civility for the achievements of others. To a specific question, a number is sent back, other participants can check if it fits the question, or send another number, theoretically without any discussion.

This form of interaction is very well known by the historians. It has been described either as a rather infantile and derivative feature of these amateurs or (sometimes and) as a trace of the necessary competition between mathematicians inside the patronage-system. Several metaphors indeed are used by Descartes in the letters concerning these challenges: game, trial, fence-fight, all related to the various situations of civil oppositions in the culture of the persons involved. However, the challenges do not touch, as participants or as judges, the patron-part of the network. Non-satisfying answers are often politely accepted, as seen for those of Descartes to Sainte-Croix' challenges, or if not, precisions or clarifications simply required. In more than a few cases, challenges degenerate into refined exercises suggested for the pleasure of the other, like between Fermat and Frenicle around the forties.

I would suggest here that challenges play specific roles and fit well with various constraints and features of the network. It helps first to solve the need for recognition achieved in communicating one's own results and the difficulty of trust in the network, specially vis-à-vis the Parisian milieu: Fermat for example protests to be as "quelqu'un de ceux du lieu où vous êtes qui s'attribue impunément les inventions d'autrui après qu'elles lui ont été communiquées,"[28] but the accusation seems to have been a current one, especially easy to justify because very few things were published and disseminated outside the very space of the correspondence. Then, it provides a suitable stimulation to work. One can contrast in this respect Descartes' vague answer in 1631 about double numbers as 120 and the acceptance to work, as seen above, when the same range of questions comes in a challenge with M. de Sainte-Croix and Frenicle. Last, but not least, the challenge is an important means of obtaining a true estimate of the difficulty of a question, to fix a value for it, as well as value for the author. Fermat to Roberval: "J'ai trouve un très grand nombre de belles propositions. Je vous envoierai la démonstration de celles que vous voudrez: permettez-moi néanmoins de vous prier de les essayer plutôt et de m'en donner votre jugement."[29] This problem—to establish the value of a result or a problem—is a delicate one, especially because of the variety of motivations and sources we have already noticed about the participants. No external system or at least a specific subsystem (as for instance referees) is delegated to judging much of the matters discussed in the correspondence—while some of the participants were official experts for issues raised elsewhere, for instance the question of longitudes.

[28] Fermat, *op. cit. in* n.9, II, 207.
[29] Fermat, *op. cit. in* n.9, II, 74.

Transitions to challenges to open questions can occur when personal trust has been established, but we have too few cases to be able to say more on the precise mechanisms: Fermat and Frenicle, possibly also Saint-Martin, enter for some times in a combination of the two exchanges. But the specificity of the challenge in the procedure of testing is underlined by several participants. When it is suggested to test Descartes' method through open questions, he answered back: "Ce n'est pas le style des géomètres que de poser des questions qu'ils ne peuvent résoudre."

Transmission of the method (or at least of a rule, of a means to find the numbers) is a worthy possession, not to be given without discernment, but to be first made valuable For instance, Fermat to Carcavi:

"Il n'y a certainement pas quoi que ce soit dans toutes les Mathématiques plus difficile que ceci et hors M. de Frenicle et peut-estre Mr Descartes, je doute que personne en connaisse le secret, qui pourtant ne le sera pas pour vous non plus que mille autre inventions dont je pourrai vous entretenir une autre fois."[30]

Challenges are precisely the form to test the method without giving it: they allow to orchestrate a comparison of them directly according to their power to find solutions (in numerical form, which can be checked and easily exchanged if necessary) to the problems of other members and to propose to them problems they cannot solve.

Allied with the emphasis on problems, it nudges mathematics in given directions:[31] to imagine complicated, but effective problems, with numerical solutions, preferably rare ones, with a great number of digits. The bigger a possible solution should be, the most efficiently it displays the knowledge and shows that it was not obtained by a pure trial-and-error approach. Within these limits, we assist to a certain kind of progress: if around 1631, the discussion on double numbers does not seem to provide much feedback, by 1634, various examples are known, by 1638, systematic processes to find all kind of submultiples of their sums of divisors are claimed by several persons, who exhibit as proofs lists of huge examples. In the meanwhile, Euclid's procedure to find perfect numbers—numbers equal to the sum of their divisors—has been scrutinized, new examples discovered, shortenings explored. These questions lead Fermat and Frenicle in particular to new explorations, which we are tempted to describe as deeper and more theoretical, for instance the identification of the form of divisors of certain type of expressions, as a power of 2 plus or minus 1, and even (false) conjectures. From 1642 on, some of these new questions about prime numbers circulate around, while the newcomers in the area, like de Villiers or Thibaut, still ask Mersenne for the computation of the submultiples, the sum of the divisors of a numbers, etc. That is, we have a clear pattern of difficult questions solved (at least partially) through a collective effort, giving raise to new, more crucial problems, and whose solutions were disseminated to the whole circle—indeed even to a more general public, through Mersenne's books.

[30] Mersenne, *op. cit. in* n.5, XII, 249.

[31] Effects on individual achievements have been discussed in detail for Fermat and Frenicle, in Goldstein, *op. cit. in* n. 8 and *op. cit. in* n. 21, respectively.

To summarize, the economy of challenges efficiently structures the production of mathematical knowledge in the correspondence: it stimulates original work, it allows an evaluation of one's own and others' achievements while managing the tension between recognition and trust, it incorporates the issues of time, courteousness and integration in a heterogeneous environment. And last, but not least, it deals with the key-problem of the basic knowledge common to all, by reducing the mathematical exchanges to their common parts.

III—Challenges and controversies

This mathematical organization bears upon the issue of controversies. They appear as pulverized all over Mersenne's correspondence, both uniformized and scattered into tiny micro-conflicts. It is then quite difficult to locate what is a controversy, a simple misunderstanding or a priority dispute.

How revealing for our topic is for instance the rebuff that Fermat meets with when trying to answer a challenge of André Jumeau de Sainte-Croix in 1636: "Trouver deux nombres, chacun desquels comme aussi la somme de leur agrégat ne conste que de trois tétragones?"[32] Fermat proposes to take two squares the sum of which is also a square (for instance 9 and 16) and to multiply each of them by a number sum of three squares, like 11, and adds immediately how to find an infinity of such solutions, as well as generalizations. According to Mersenne, the answer is not satisfying because Sainte-Croix wanted numbers which could not be decomposed also into a sum of four squares. Fermat counters thanks to the double authority of a plural and of a classical reference:

"quand nous parlons d'un nombre composé de trois quarrés seulement, nous entendons un nombre qui n'est ni quarré, ni composé de deux quarrés; et c'est ainsi que Diophante et tous ses interprètes l'entendent."[33]

Sainte-Croix' question in his sense would be solved by chance, not by "une conduite assurée" (that is, a rule, or a method). The formatting of the challenges hampers us to understand what is here in question: a terminological misunderstanding, a difference of cultures and sources, the wish to test some deeper issue?

[32] Fermat, *op. cit. in* n.9, II, 49.
[33] Fermat, *op. cit. in* n.9, II, 77.

The situation does not seem very different a few years later when Frenicle challenged Mersenne's network, and in the first place Descartes, with the "question des ellipses," as it will soon be called:

"On demande qu'une ligne rationnelle AB serve de grand diamètre a tant d'ellipses qu'on voudra (& non plus), dont la surface convexe du verre [HBO] dont la face qui est au dehors, HO, soit circulaire dont le centre F' soit le point brulant exterieur de l'ellipse (lequel centre doit aussy etre la place de l'objet) sera portion, auxquelles le petit diametre CD soit moindre que la distance des points brulants FF', & que tant le petit diametre de l'ellipse, que l'epaisseur du verre OB, & la distance de l'objet a iceluy F'O soient lignes rationnelles."[34]

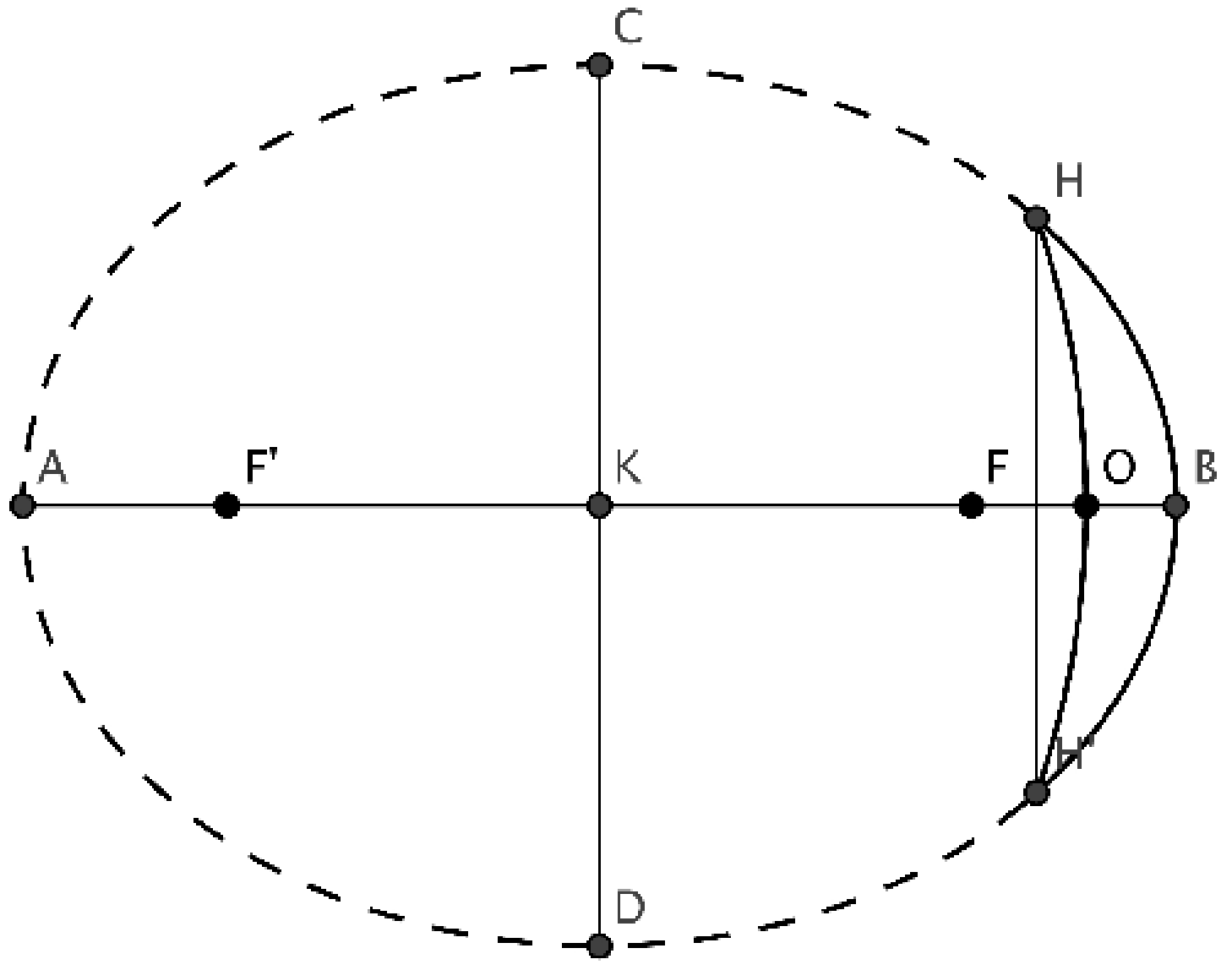

Despite a superficial appearance of an optical challenge, the problem was arithmetical: all the sought-for quantities are rational numbers, and in fact integers. The disguise had a double function: it was both an act of courtesy towards correspondents who were not usually working on numbers (a fact perfectly recognized by Descartes, for instance) and a complexification to discourage too ignorant people. Reactions to the challenge varied nonetheless from annoyance to several kinds of rewritings transforming it into a problem on Pythagorean triangles or into an equation and solving it partially. Descartes, in particular, used of course his algebraic notations to reformulate the problem: choosing OF =a, BF=b as his two main unknowns, he obtained AB= 2bb/a et CD=2b√(2b/a)-1 and sent this as a "general solution" to Frenicle. It was rejected, and it may appear at a superficial glance that the reason is of the same sort that Sainte-Croix used to disqualify Fermat's ones in 1636: Frenicle insisted on the supplementary condition, "que cette ligne ne servist point a plus grand nombre d'ellypses qu'a celuy qui seroit demandé." Descartes' solution indeed did not identify precisely when the AB and CD are rational numbers, and a fortiori how many ellipses correspond to AB. Indeed, from an arithmetical and combinatorial point of view, this was where the difficulty lied. For Descartes, it was not the case:

[34] On this problem, see Goldstein, *op. cit. in* n. 21.

"Ayant trouve tout d'abord tout ce qui me semblait contenir de difficulté en la question, qui estoit de donner autant d'ellipses rationelles qu'on voudrait qui eussent une mesme ligne pour plus grand diamètre, et ayant d'autres pensées dans l'esprit, je ne me suis pas arresté a considerer toutes les exceptions qu'il falloit faire, afin que cette ligne ne servist point a plus grand nombre d'ellypses qu'a celuy qui seroit demandé."[35]

For Descartes, the priority is the universality of algebraic relations, enumerative issues are thus handled as "exceptions" (thus dismissed). For good reasons, Frenicle did not accept this interpretation, and the exchange became acrimonious. Descartes ended it by disqualifying both the problem and the author:

"Je n'avais plus envie de répondre, car sa question n'est ni belle ni industrieuse. […] mais je ne veux point contester, car il me paroist estre […] du nombre de ceux qui veulent, a quelque prix que ce soit, avoir gaigne & parler les derniers."[36]

The correspondence does not offer much more, but if we follow the protagonists in their other texts, elements of a quite deep controversy emerge: the multifarious one which followed the publication of Descartes' Géométrie, and touches upon the personal, the disciplinary and the cognitive levels at the same time. Frenicle, most particularly, developed challenges that questioned the supremacy of symbolic, Cartesian algebra within the field of arithmetic.

A useful hint to capture such issues is to look for cases where the communication is broken: one correspondent announces, exactly as Descartes did, that he will stop answering the challenges. It is a remarkable feature of the correspondence that such statements are not rare: people address their complaints to a third person, often Mersenne, who in turn informs the guilty part of the pending situation. Such a disruption, indeed, is a rather drastic move, as it threatens the very existence of the correspondence as a mathematical institution: in Mersenne's network, silence has to be negotiated, explained, compensated. Something, anything, has to be written back, and it is a priority with respect to, for instance, mathematical requirements. "J'aime mieux paraitre ignorant en vous répondant mal qu'indiscret en ne vous répondant point du tout"[37], wrote Fermat to Mersenne.

Following these clues, we found two main situations. One is at work in the Frenicle-Descartes' controversy on the "question des ellipses", it touches upon a disciplinary issue. Problems on numbers, arithmetic in itself are often disregarded, as fastidious, without use, and requiring no subtlety. The blind use of tables, or the possibility to rely on simple trials, one number after another, to reach the solution, without method, are current criticisms against most problems on numbers. On the question of ellipses, de Beaune comments:

---

[35] Descartes, *op. cit. in* n. 15, II, 536.
[36] Descartes, *op. cit. in* n.15, II, 567.
[37] Fermat, *op. cit. in* n.9, II, 3.

"Je vous supplie de me dispenser de la recherche de ceste question, pour m'apliquer, aux heures de mon loisir, a de plus serieuses: cette question n'estant d'aulcun usage et ne tombant poinct soubs la science des rapports, qui les considere universelement, aussi bien entre les lignes commensurables et incommensurables, si bien que la recherche en seroit extremement laborieuse et de nul proffict, ce qui n'arrive pas en celles de geometrie et celles d'arithmetique qui tombent soubs ceste science des proportions, les aultres estant de peu de consideration et n'estant d'aulcun usage."[38]

Or, in another context, Deschamps writes:

"Je ne doubte point cependant que sans les conter [les carrés magiques] l'on n'en puisse determiner le nombre; mais s'il n'a plus d'utilité que le seul plaisir de contenter la curiosité, le temps me sembleroit mal employé en ceste recherche."[39]

To the "groping" arithmeticians, the analysts using algebra oppose the universal power of their favorite tool. Reciprocally, the arithmeticians shape their problems in order to make them intractable by an algebraic approach. "Je sais que l'algèbre de ce pays n'est pas propre pour soudre ces questions," comments Frenicle.

The second situation in which a rupture is announced tackles epistemological issues. In 1643, Fermat for instance, sent a challenge to Frenicle and Saint-Martin, in particular: to find a right-angled triangle such that the longest side (the hypotenuse) is a square [of an integer] and the sum of the two smallest ones is also a square [of an integer]. Shortly after, Fermat learnt, from Mersenne, that "[s]es questions impossibles ont refroidi M. Frenicle et Saint-Martin." Since the beginning of his participation, Fermat indeed has challenged others with impossible questions—for instance, that no rational cubes can be the sum of two rational cubes, a particular case of the celebrated so-called Fermat Last Theorem. Such problems have not the adequate form to circulate smoothly in the network: they have no numerical solution that can be exhibited, they require a proof. Fermat tried by several mediators to convince Frenicle to come back into the exchange. But in order to succeed, he had finally to display a numerical solution, that is to show that he had not betrayed the trust of his correspondents by an unfit challenge: the smallest solution is the triangle with the sides a= 1 061 652 293 520, b= 4 565 486 027 761, c=4 687 298 610 289… If in this case, the challenge did not degenerate into a conflict, the conditions in which it could have been are clearly visible: the permissible statements are those that fit the conditions of the correspondence. Fermat would even work hard, so he said, in order to be able to apply his favorite method of proof (by infinite descent, a proof by absurdum a priori better adapted to show that some relation among integers is impossible) to positive statements.

Disciplinary and epistemological issues are of course tightly intertwined. In his short-lived controversy with Frenicle, there is a misunderstanding on the nature of which of his methods

[38] Mersenne, *op. cit. in* n.5, VIII, 360.
[39] Mersenne, *op. cit. in* n.5, VIII, 545.

Fermat is trying to test: Frenicle and Saint-Martin believed that it was again the proof of impossible problems, while Fermat was in fact testing here the possibility of algebra in such numerical questions: his discovery of the triangle comes from a subtle use of algebraic transformations inside Diophantine problems. [40] Reciprocally, the Descartes-Frenicle controversy displays a type of statement which is both adapted to the constraints of the correspondence as we have seen them, and apt to put on trial Cartesian algebra: the question of exact enumeration. "Trouver une ligne […] qui serve de diametre à tant d'ellipses et non plus," or as Saint-Martin requires from Fermat at another occasion: "déterminer combien de fois 1 803 601 800 est la différence des (plus grands) côtés d'un triangle qui ait un quarré pour différence de son petit côté aux deux autres côtés."[41] Questions involving "exactly how many" became a trade-mark of the arithmeticians participating to Mersenne's correspondence network: they allowed the participants to discriminate those who possessed a serious knowledge of the properties of integers, while respecting the rules of the communication.

## IV—Routine controversies

Mersenne's correspondence, as a mathematical institution, fulfilled its main objective: it socialized together persons interested in various degrees by mathematics, but coming from extremely different social circles and scientific backgrounds. Mathematics were used a social link, of course, but the model is not one of strict separation of roles and tasks. Agreement came first of all about a largely shared organization of mathematical practice—or even more, of its image: effective (interpreted as practical) solutions and rules of conduct were both at the top of the priorities, in a configuration reminiscent of Baconian ideals. To test methods—and even more one's own method—in the conditions of restricted trust that pervaded the network required to go through the solution of problems, that in turn were the source of evaluation for any mathematical results.

The constitution of that institution did not rely on the establishment of a canon, for instance, which would identify legitimate sources and important problems, it was reinforced bits by bits through each concrete numerical answer. As far as arithmetic is concerned, the opposition to algebra, by a large part of the Parisian mathematicians, favored the development towards certain directions, to the cost of other aspects, in particular the administration of proofs.

These features, both social and intellectual, appear as archeological components of mathematical practice. Challenges, not controversies as such, stimulated the work to do. Controversies, indeed, were trapped into the concrete shapes of this collective structure. Most of them appear thus, paradoxically, as routine controversies, operating inside normal science. If some of them—in particular that born from the resilience to symbolic analysis—could probably have a larger effect in the long-term, their potential effect was muffled inside the correspondence in the simple exchange of problems and numbers.

[40] See for details Goldstein, *op. cit. in* n.8.
[41] Fermat, *op. cit. in* n.9, 250.

Catherine Goldstein
Directrice de recherche au CNRS
Histoire des sciences mathématiques
Institut de mathématiques de Jussieu – UMR 7586
Case 247
4, place Jussieu
F-75252 Paris Cedex, France.
cgolds@math.jussieu.fr